\documentclass[12pt]{article}
\begin{document}
\def\pni{\par \noindent}
\def\vsh{\smallskip}
\def\vs{\medskip}
\def\vvs{\bigskip}
\def\vvvs{\bigskip\medskip} 
\def\vsp{\vsh\pni}
\def\vsn{\vsh\pni}
\def\cen{\centerline}
\def\ra{\item{a)\ }} \def\rb{\item{b)\ }}   \def\rc{\item{c)\ }}

\sloppy

\begin{center}

{\large\bf George William Scott Blair -- the pioneer of \\
factional calculus in rheology}\footnote{Short version published at 
{\it Communications in Applied and Industrial Mathematics}. -- 2014, v. 6, no. 1. -- e481. {\tt ISSN 2038-0909; DOI: 10.1685/journal.caim.481} }

\vspace{3mm}
{\bf Sergei Rogosin}$^{1, 2}$, {\bf Francesco Mainardi}$^3$

\vspace{3mm}
{$^1$IMPACS: Institute of Mathematics, Physics and Computer Sciences,  
Aberystwyth University, Wales, UK\\
ser14@aber.ac.uk}

{$^2$Department of Economics, 
Belarusian State University, Minsk, Belarus \\
rogosin@bsu.by}

{$^3$Department of Physics, University of Bologna,  and INFN, Italy \\ 
francesco.mainardi@bo.infn.it}

\end{center}

\begin{abstract}

The article shows the pioneering role of the British scientist, Professor G.W.Scott Blair,
in the creation of the application of fractional modelling in rheology. Discussion of his results is presented.
His approach is highly recognized by the rheological society and is adopted and generalized by his successors.
Further development of this branch of Science is briefly described in this article too.

\end{abstract}


\section*{Foreword to the second extended version.} 

After appearance of the first version of the paper in arXiv it has got an interest as from Mathematical as from Rheological Society. 
In March 2014 one of the authors (SR) made a report at the research seminar at Aberystwyth University. 
The content of this seminar was described later by Prof. Simon Cox in his note published at Rheological Bulletin \cite{Cox}. 
For readers' convenience we present this note in the Appendix to the extended version of this ArXiv paper. We also add some additional references justifying the role of Scott Blair's results (\cite{Jai14}, \cite{JaiMcK13}, \cite{MaiSpa11}). 
On this respect it is interesting to mention the ``Bingham Lecture'' \cite{McKJai} by G.H. McKinley. We pay attention of the readers to our recent book \cite{GKMR} devoted to the description of the properties of the Mittag-Leffler function, the Queen  Function of the Fractional Calculus.  

\section{Introduction.}

The interest to applications of fractional calculus in modelling of different phenomena in Physics, Chemistry, Biology is rapidly increasing in the recent three decades. First of all we have to point out the constitutive modelling of non-Newtonian fluids. The main reason is that the fractional models give us possibility for simple description of complex behaviour of a viscoelastic material. In  pioneering (mainly experimental) works of 1940-1950th, it was
discovered, for instance, that the relaxation processes in some materials exhibit an algebraic decay, which cannot be described in the framework of the Maxwell model based on exponential behaviour of the relaxation moduli.

In order to see perspective in the development of fractional models it is important to understand how such models appear, and what was really done by the pioneers. Among the works  which are in the core
of the first period of the fractional modelling one can single out the series of articles and monographs by G.W. Scott Blair.
His role is not overestimed by the rheological society (see, e.g.  \cite{Bar99}, \cite{Dor02}), but anyway, some details of his work are still of great importance. We propose here an analysis of the results by G.W. Scott Blair along with their influence on the modern development of the fractional modelling in rheology.

\subsection{Short biography by G.W. Scott Blair} Dr. George William Scott Blair (1902--1987) was born on 23 July 1902, of Scottish ancestry, in Weybridgein Surrey, England. After graduated a famous public school at Charterhouse he went to Trinity College Oxford in 1920, where he studied Chemistry, with Prof. Sir Cyril Hinshelwood as his tutor. He carried out a one-year research project in colloid chemistry to complete his master thesis with honour degree. 

 After graduating Scott Blair was employed as a colloid chemist with a Manchester firm of Henry Simon, working there on the viscometry of flour suspensions, publishing his first rheology paper in 1927. In 1926 he was offered a post in the Physics Department of the Rothamsted Experimental Station, where he was working on the flow properties of soils and clays until 1937.
 It was there, where he made with his colleagues the first quantitative study of so called sigma-phenomenon, which was originally described by Bingham and Green in 1919. Schofield and Scott Blair (see \cite{SchSco30}) studied this phenomenon from 1930 at Rothamsted for soil and clay pastes and named it ``sigma effect''. These studies were probably unknown to F${\dot{\textrm a}}$hraeus and Lidquist, who first discovered the sigma effect for blood, referred to as the ``F${\dot{\textrm a}}$hraeus-Lindquist phenomenon''.
 At this period some preliminary experiments were provided by Scott Blair which led him later to the necessity to consider anomalous relationship between stress, strain and time (see, e.g., \cite{SchSco33}).

In 1929 Scott Blair took up a Rockefeller Fellowship at Cornell University in Ithaca, New York state. He worked there on the flow of potter's clay and developed a means of measuring its plasticity. He attended the inaugural meeting of Society of Rheology in December 1929 in Washington DC, met there all pioneers of Rheology, Eugune Bingham and Markus Reiner among them, and many became his life-time friends.

After returning to Rothamsted, Scott Blair made rheological research on honey and flour doughs. He also studied  together with the well-known psychologist David Katz, psychophysical problems in bread making. His interest in psychology led him, together with F.M.V. Coppen, to initiate a new field, for which he coined the word ``psychorheology''. It is considered as one of the fields of biorheology.

In 1936 he submitted his PhD thesis to the University of London, and it was examined by Prof. Freundlich. The same institution later awarded him a D.Sc. for his labours in rheology (probably the first ever rheology D.Sc.) In 1937 he joined National Institute for Research in Dairying, University of Reading as a head of Chemistry but soon took over the newly formed Physics Department and remained in that position until his retirement thirty years later.

In 1940 the British Rheological Society was founded.
Scott Blair played a prominent role and took active part in the development of rheology. He was a Founder-Member and, later, president of the British Society of Rheology (1949--1951). He took part in the organization of the First International Congress on Rheology, held at Scheveningen, Holland, in 1948. He was a Secretary of the Second International Congress on Rheology in Oxford in 1953 and a member of Committee on Rheology, set up by the International Council of Scientific Union.
 Scott Blair was given to flights of fancy into psychorheology, fractional differentiation etc. In 1969 he was awarded the Poiseuille Gold Medal of the International Society of Haemorheology (now Biorheology) and in 1970 he received the Founders Gold Medal of the British Society of Rheology.  Together with J. Burgers he published a monograph on rheological nomenclature \cite{BurScoB49}. For many years Scott Blair was the Chairman of the British Standard Institute Committee on Rheological Nomenclature.

During almost a half of  a century, George W. Scott Blair was  one of the leading rheologists. Beginning from 1957 Scott Blair devoted his experimental and theoretical work entirely to hemorheology. Since the foundation of the International Society of Hemorheology in Reykjavik, Iceland in 1966, he was a member of its Council and acted as Chairman of its Committee on Standards and Terminology. After he retired he worked on the flow and coagulation of blood at the Oxford Haemophilia Centre.

Scott Blair was very active in publication and editorial work. He  was a co-founder of the Journal ``Biorheology'' and its Co-Editor-in-Chief from its inception in November 1959 to December 1978 (see \cite{Cop88}). The books and research papers of Scott Blair were donated to the British Society of Rheology and later deposited in the Library of Aberystwyth University in early 1980's. The collection has over 550 books and its aim is to develop this ``into an up-to-date library of rheological literature available to all members of Society''. Rheology Abstracts and the British Society of Rheology Bulletin are two journals published by/for the Society which form an important of the Collection. The books and journals catalogued online (access via http://primo.aber.ac.uk).

\section{Rheology and Psychophysics}

It was Professor Bingham who had chosen the name ``Rheology'' for this branch of the Science and
gave the definition of it: ``The Science of Deformation and Flow of Matter'' (see \cite{ScoB47a}) motivated by Heraclitus' quote
``$\pi\alpha\nu\tau\alpha\; \rho\varepsilon\iota$'' (``everything flows''). Rheology is one of the very few disciplines
having exact day of its birth, April 29, 1929, when the preliminary scope of the Society of Rheology was set up by a
committee met at Columbus, Ohio. Anyway, the ancient Egyptian scientist Amenemhet (ca 1600 BC), who made the earliest application
of the viscosity effect,  can be considered as the first rheologist (see, {\it e.g.}, \cite{ScoB44a}).

The observables in rheology are deformations or strains, and the changes of strains in time. Changes of strains in time
constitute a flow. Thus, these changes are generally associated with internal flow of certain kind. States of stress  are
inferred either from the comparative strain behaviour of complex and simple systems in interaction or from the behaviour of a
known mass in the gravitational field.  In physical testing, stresses $(S)$, strains  $(\sigma)$ or their differentials with respect to time $\left(\dot{S}, \dot{\sigma}\right)$ are normally held constant, leaving either a length to be measured, or the time $(t)$.

There is a group of fluids which is characterized by a coefficient of viscosity for a specific temperature. These fluids, known as Newtonian fluids,
were singled out by Newton who proposed the definition of resistance (or viscosity in modern language) of an ideal fluid. Pioneering
work on the laws of motion for real (i.e. non-ideal) fluids with finite viscosities was carried out by Navier \cite{Nav23} and later by Stokes \cite{Sto66}.
The Navier-Stokes equation enabled, among other things, prediction of velocity distributions and flow between rotating cylinders and cylindrical tubes (see \cite{Dor02}).

Nowadays rheology generally accounts for the behaviour of non-Newtonian fluids, by characterizing the minimum number of functions that are needed
to relate stresses with rate of change of strains or strain rates. This kind of fluids is called Newtonian since
 Newton's introduction the concept of viscosity.

In practice, rheology is concerned with extending continuum mechanics to characterize flow of materials,
that exhibits a combination of elastic, viscous and plastic behaviour by properly combining elasticity and (Newtonian) fluid mechanics.
In \cite{Dor02} the main directions in the development of the rheology are described. First of all this is linear {\it viscoelasticity}. One of the most important early
contribution in this area is the work by Maxwell \cite{Max68}. In order to explain the behaviour of the materials which are neither truly elastic nor viscous
he proposed a constant relaxation time $(t_r)$ and justified implicitly the model of a dash-pot and spring in series. Anyway, he realized that for some materials
the assumption of constancy of  the relaxation time is over-simplification, in these cases $t_r$ has to be a function of stress. Meanwhile, the notion of Maxwell's
units (i.e. pieces of a material having constant relaxation time) has been widely explored by rheologists.
Later the conception
of the ``orientation times'' $\tau$ has been developed (see, {\it e.g.}, \cite{AleLaz40}). It is considered unit stress conditions, supposing that the strain is approaching to an equilibrium value. Thus, the immediate Hookean strain is first subtracted and  $\tau$ is defined as the time taken for the remaining strain, resulting from the orientation of the chains.

 Another direction which was singled out in \cite{Dor02} is the study of {\it generalized Newtonian materials}.  This type of fluid behaviour is associated with the work by Bingham \cite{Bin22}
 who proposed so called yield stress to describe the flow of paints. In \cite{ScoB-Cop42a}, it  has been pointed out the close similarity between the usual experimental Bingham curve and the curve of a high power-law. Thus, it shows possibility of existence of systems for which the Bingham plot gives a fairer and simpler account of the data.

 The study of {\it non-linear viscoelasticity} started at the begining of XXs century, when the area of rheology was most rapidly grown (see, {\it e.g.} \cite{Dor02}).
 Thus, Poynting \cite{Poy13} in his experiment discovered that loaded wires increased by a length that was proportional to the square of the twist, what did not correspond
to the usual expectation of the linear viscoelasticity theory. Probably the first theoretical work on non-linear viscoelasticity was done by Zaremba \cite{Zar03}, who extended linear theory to the non-linear regime by introducing corotational derivative in order to incorporate a frame of reference that was translating and rotating with the material. More extended description
 of the results in non-linear viscoelasticity can be found in \cite{Dor02} (see also \cite{Mai10}, \cite{Uch13b} and references therein).

Not all properties of flowing matter can be interpreted in term of real rheological sense. In this case psychophysical approach with its
psychophysical experiments can be helpful. Psychophysics is defined as the scientific study of the relation between stimulus and sensation (see, {\it e.g.} \cite{Ges97}).
Psychophysicists usually employ experimental stimuli that can be objectively measured. Psychophysical experiments have traditionally used three methods for testing subjects' perception in stimulus detection and difference detection experiments: the method of limits, the method of constant stimuli and the method of adjustment.

G.W. Scott Blair widely used psychophysical experiments in his research (see \cite{ScoB74}). Therefore, it is interesting to recall how he described the role of psychophysics in rheology (see \cite{ScoB47a}): ``The complex and commercially important rheological ``properties''
 of many industrial materials are still assessed subjectively by handling in factory and are expressed in terms ``body'', ``firmness'', ``spring'', ``deadness'', shortness'', ``nerve'', etc. - concepts which cannot be interpreted ... in terms of simple rheological properties at all. In view of this fact ... it is clearly advisable to know something of the accuracy with which these handling judgements can be made and, by bulking sufficiently large numbers of data together so that reproducibility is ensured, to attempt to correlate the entities so derived with manageable functions of $S : \sigma : t$. A start has been made in this direction and not only have a number of reproducible regularities been observed, but the information obtained has laid the foundations of a theory of
 ``{\it Quasi-properties}'' which it is hoped will facilitate the study of purely ``physical'' rheology of complex materials.''

 This observation is a core of Scott Blair's method which he used along his career.

\section{Scott Blair Fractional Element}

\subsection{Nutting's Law}

In 1921 Nutting reported (\cite{Nut21}) about his observation that
mechanical strains appeared at the deformation of the viscoelastic materials
decreasing as power-type functions in time. From a series of experiments, which covered a range of materials from the elastic solid to the viscous fluid, Nutting suggested a general formula relating shear stress, shear strain and time, whenever shear stress remains constant:
\begin{equation}
\label{Nutting}
\sigma(t) \sim C \Delta S \cdot t^{-\alpha},
\end{equation}
with constant order $\alpha\in (0, 1)$ which is close to $1/2$ for many materials.
This conclusion was in a strong contradiction to the standard exponential law.
Later the Nutting's observation was justified by Gemant who studied the properties of
viscoelastic materials under harmonic load. It was shown that the memory function $\eta(t)$
can have power-type relaxation behaviour proportional to $t^{-3/2}$. In 1950 Gemant published a series of 16 articles entitled ``Frictional Phenomena'' in
Journal of Applied Physics since 1941 to 1943, which were collected in a book of the same title
\cite{Gem50}. 
In his eighth chapter-paper 
\cite[p. 220]{Gem42}, he referred to his previous articles \cite{Gem36}, \cite{Gem38}
for justifying the necessity of fractional differential operators to
compute the shape of relaxation curves for some elasto-viscous fluids.

Gemant has used half-differential, but is his later papers he says that
fractional differential ``only occurs as a useful mathematical symbol, whereas
the  underlying elementary process, whatever it may be, will probably contain
differential quotient of an integral order''.

Scott Blair surely knew the attempts by Gemant (see, {\it e.g.} \cite{ScoB47b}) to generalize
Maxwell's theory by changing various (but integer) powers in complex modulus of the Maxwell
Fluid Model to fractional powers. In fact, Scott Blair (together with Coppen) also came to the form of Nutting equation, but from another consideration. They argued that the material properties are determined by various states between an elastic solid and a viscous fluid, rather than a combination of an elastic and a viscous element as proposed by Maxwell.
In \cite{ScoB-Cop39} it was  pointed out that, since for Hookian solids, strain is proportional to stress and to
unit power of time, for intermediate materials, it might be expected to be proportional to stress and to
some fractional power of time with exponent $\alpha, 0 < \alpha < 1$ and described this relation in the form
\begin{equation}
\label{Nutting-SB}
\psi = S^{\beta} \sigma^{-1} t^{\alpha},
\end{equation}
where proportionality coefficient $\psi$ is a constant. Derived in this way the equation (\ref{Nutting-SB}) looks
entirely empirical, tough the fundamental significance of $\alpha$ (which is called the {\it dissipation coefficient})
is shown in psychophysical experiments described by Scott Blair and Coppen (see \cite{ScoB-Cop40}, \cite{ScoB-Cop42}, \cite{ScoB-Cop43}).

A comparison of the partially differentiated Nutting equation and Maxwell's equation may be written  (see \cite{ScoB44}), namely, for Nutting:\footnote{The suffix $\sigma$ indicates shear strains.}
\begin{equation}
\label{Nutting_dif}
- \left(\frac{\partial S}{\partial t}\right)_{\sigma} = \frac{\alpha}{\beta} \frac{S}{t},
\end{equation}
and for Maxwell:
\begin{equation}
\label{Maxwell_dif}
- \left(\frac{\partial S}{\partial t}\right)_{\sigma} =  \frac{S}{t_r}.
\end{equation}
Since Nutting equation gives $t = \psi^{1/\alpha} S^{-\beta/\alpha} \sigma^{1/\alpha}$, it is apparent that the Nutting
treatment postulates a single {\it relaxation time} proportional to a power of the stress. This is simplest possible way of implementing Maxwell's suggestion that
relaxation time $t_r$ may be some function of stress. From the other side it justifies the believe that the use of fractional calculus in description of processes toward equilibrium  is necessary if one has to keep the Newtonian time scale.

\subsection{Scott Blair's fractional model}

It was suggested in \cite{ScoB-Cop43} that, in considered cases, comparative firmness
is judged neither by $\sigma$, nor by $\dot{\sigma}$, nor by any mixture of these two
magnitude, but by some intermediate entity, namely by fractional derivative
 $\frac{\partial^{\nu} \sigma}{\partial t^{\nu}}$.\footnote{In fact, few misprints had been made in  \cite{ScoB-Cop43}
later corrected by Scott Blair. Thus, in the original reprint of \cite{ScoB-Cop43} in the Scott Blair reading room one can find hand-written corrections made by Scott Blair.} More exactly, he wrote: The general constitutive equation ``... is applicable to integral values of $n$ but a more general equation may be used even $n$ is a fraction.The numerical coefficient is expressed as a quotient of $\Gamma$-functions and may be written
\begin{equation}
\label{ScoB-Cop43}
\frac{\delta^n \sigma}{\delta t^n} = \frac{\Gamma(k + 1)}{\Gamma(k - n + 1)} t^{k - n} \Psi^{-1} S.
\end{equation}
The expression $\Gamma(k + 1)$  is given by
$$
\int\limits_{0}^{\infty} e^{-x} x^{k} dx
$$
and, whatever the value of $n$, the  $\Gamma$ factor is a number independent of $\Psi_{m}$ and $t_c$, so that the validity of the plot is unchallenged.''

In his work Scott Blair did not specify what kind of fractional derivative he used. From the way how he has calculated derivative of any power we can conclude that this is the standard Riemann-Liouville derivative.
It is   quite instructive to cite some words by Scott-Blair  quoted by Stiassnie in
their correspondence, see  \cite{Stia79}:
  {\it I was working on  the assessing of firmness of various materials
 (e.g. cheese  and clay by experts handling them) these systems are of course both elastic
 and viscous but I felt sure that judgements were made not on an addition of elastic and viscous
 parts but on something in between the two so I introduced  fractional differentials of strain
 with respect to time}.
 Later, in the same letter  Scott-Blair added:
{\it I gave up  the work eventually, mainly because I could not find a  definition of a fractional
differential that would satisfy the mathematicians.}

The above said  {\it Principle of Intermediacy}
 was discussed in details by Scott Blair in \cite{ScoB44} basing on purely physical grounds.
The theory of fractional modelling in rheology is developed by Scott Blair, Veinoglou and Caffyn in \cite{ScoB47b}.
 In \cite[p. 30]{ScoB47a} it is briefly summarized: ``... times are normally defined as equal when ``free''
 Newtonian bodies (or alternatively light) traverse equal (superposable) distances in them. This leads to
 a a definition of velocity as the first differential of length with respect to time which,
 because of this definition of time equality, is constant for Newtonian bodies; and to the
 second  differential, called acceleration.

 When bodies are not influenced by other bodies, and their velocities change with time, a force is postulated
 and defined as rate of change of (velocity $\times$ mass). It is long been realized that the Newtonian time scale arbitrary (see
 \cite[p. 80]{Poi04}) and in the case of a complex plastic being strained, the rheologically active units are certainly
 not independent Newtonian bodies. It should, therefore, be easy to choose a non-Newtonian time
 equality definition which would reduce the entities by which firmness is judged to simple
 whole-number\footnote{I.e. non-integer order.} differential expression. The use of separate time scales for different materials is not convenient, however, so Newtonian time is used, but, as a result of this arbitrary procedure, the derived constants cannot be expected to be built up entirely from whole-number differentials. It is thus apparent that fractional differential is an essential feature of our whole mode of approach.''

 In \cite{ScoB-Caf49} are discussed the circumstances under which it is practicable to express the Nutting equation and its fractional derivatives
 in a simple dimensional form. Three main principles of a new proposal are formulated: (1) the fact that the treatment does not lead to any understanding
 of structure of the materials or of their molecular configurations; (2) the only entities are used whose dimensions depend of the nature of of the material; (3) fractional derivatives and corresponding coefficients are understood as something intermediate between zero and first derivatives and corresponding coefficients.
Scott Blair highly supported (see, {it e.g.} \cite{ScoB-Cop42a}) the ideas by Nutting supposing that for that moment it describe a special but very frequently adequate cases. Anyway, he though that the phenomenon dealing with Nutting equation are related to the fundamental structure of materials.

 Fractional derivative of order $\mu, 0 < \mu < 1,$ with respect to time $t$  of the Nutting equation (in the form (\ref{Nutting-SB}) gives (see \cite{ScoB-Cop43}, \cite{ScoB47b})
 \begin{equation}
 \label{Nutting-frac}
 \frac{\partial^{\mu} \sigma}{\partial t^{\mu}} = \frac{\Gamma(\alpha + 1)}{\Gamma(\alpha - \mu + 1)} t^{\alpha-\mu} \psi^{-1} S.
 \end{equation}
 Another way to justify this relation is to introduce quasi-property $\chi_1$
 by the {\it Principle of Intermediacy}
 \begin{equation}
 \label{quasi-p}
 \chi_1 = S \div \frac{\partial^{\mu} \sigma}{\partial t^{\mu}},
 \end{equation}
 since the viscosity can be defined by the relation $\eta = S \div \frac{d \sigma}{d t}$, and shear modulus as $n = S \div \sigma$.

 The relation (\ref{quasi-p}) can be integrated to give
 \begin{equation}
 \label{Nutting-gen}
 \sigma = S^\beta \left(A t^{\alpha^{\prime}} + B t^{\alpha^{\prime} - 1} + C t^{\alpha^{\prime} - 2} + \ldots\right),
 \end{equation}
 where $A$, $B$, $C$ etc. are constants. Clearly this equation coincides with Nutting equation if $A \gg B, C \ldots$.

\section{Fractional models in rheology}

After the first applications of the fractional derivatives in the modelling of the processes in rheology several other fractional models were proposed to describe certain rheological phenomenon. We briefly outline here the most discussed models of such a type.\footnote{Here and in what follows we will use modern notations for stress $(\sigma)$ and strain $(\varepsilon)$ that are not be confused with the corresponding notations used by Scott Blair, namely $(S)$ and $(\sigma)$. Furthermore we will write
 $D_{0+}^{\alpha}$ to denote  the Riemann-Liouville fractional derivative implicitly  adopted by Scott Blair.}

Gerasimov \cite{Ger48} used similar arguments as Scott Blair (in, e.g. \cite{ScoB-Cop43}, \cite{ScoB47b}), namely,  interpolation between Hook and Newton's law, in order to introduce a rheological constitutive equation in terms of a  precise notation of fractional derivative
\begin{equation}
\label{Gerasimov}
\sigma(t) = \kappa_{\alpha} D_{0+, t}^{\alpha} \varepsilon(t).
\end{equation}
This equation was used in \cite{Ger48} for description of the flow of the viscoelastic between two parallel plates. He obtained an exact solution by using operational method. He has started his consideration by appealing to the Boltzmann equation saying that from experiments follow the importance of a special case of the Boltzmann equation corresponding only to the hereditary    part of the stress $\sigma(t)$
\begin{equation}
\label{Gerasimov1}
\sigma(t) = \int\limits_{0}^{t} G(t) \varepsilon(t - \tau) d \tau,
\end{equation}
or even the processes for which $\sigma(t)$ has a memory on the velocity of all earlier deformations
\begin{equation}
\label{Gerasimov2}
\sigma(t) = \int\limits_{0}^{t} K(t) \dot{\varepsilon}(t - \tau) d \tau.
\end{equation}
For the kernel in this integro-differential relation he claim that for certain materials this kernel (relaxation function) has the a form
$$
K(\tau) = \frac{A}{\tau^{\alpha}},\;\;\; A > 0, \;\;\; 0 < \alpha < 1.
$$
Hence, equation (\ref{Gerasimov2}) can be written as\footnote{Up to the constant multiplier the right-hand side
coincides with the fractional derivative known as  Caputo derivative.}
\begin{equation}
\label{Gerasimov3}
\sigma(t) = \frac{\kappa_{\alpha}}{\Gamma(1 - \alpha)} \int\limits_{0}^{t} \frac{\dot{\varepsilon}(t - \tau)}{\tau^{\alpha}} d \tau.
\end{equation}
In particular, $\alpha = 1$ gives us the Newton law, and $\alpha = 1$ corresponds to Hookean law.

Same approach is used in the above article for the study of the rotational viscoelastic flow between two concentric cylinders.

Similar to (\ref{Gerasimov}) formulation of the fractional model was proposed by Slonimsky \cite{Slo61}).

Rabotnov (see \cite{Rab48} and more extended description in his monograph \cite{Rab80}) presented a general theory of hereditary solid mechanics
using integral equations (see also \cite{Koe84}, where the use
of integral equations for viscoelasticity was revisited and interjects fractional
calculus into Rabotnov's theory by the introduction of the spring-pot was presented).

Rabotnov introduced an  hereditary elastic rheological model with constitutive equation in form of Volterra integral
equation with weakly singular kernel of special type\footnote{The integral stress - strain relationship by Rabotnov  can be re-interpreted in terms of the fractional differential constitutive equation of Zener type, see later.}
\begin{equation}
\label{Rabotnov}
\sigma(t) = E\left[\varepsilon(t) - \beta \int\limits_{0}^{t_{\alpha}} R_{\alpha}(- \beta, t_{\alpha} - \tau) {\varepsilon}(\tau) d \tau\right],
\end{equation}
where $t_{\alpha}$ is the aging time, $\alpha\in (- 1, 0]$, $\beta\not= 0$, and the kernel $R$ is represented in the form of power series
\begin{equation}
\label{Rabotnov1}
R_{\alpha}(\beta, x) = x^{\alpha} \sum\limits_{n = 0}^{\infty} \frac{\beta^{n} x^{n(\alpha + 1)}}{\Gamma((n + 1)(\alpha + 1))}.
\end{equation}
Rabotnov's kernel function $R_{\alpha}(\beta, x)$ is related to the well-known Mittag-Leffler function $E_{\alpha,\beta}(z)$ highly explored nowadays in the fractional calculus and its applications, namely
\begin{equation}
\label{Rabotnov2}
R_{\alpha}(\beta, x) = x^{\alpha} E_{\alpha + 1,\alpha + 1}(\beta x^{\alpha + 1}).
\end{equation}

Of course Scott-Blair did not know the Mittag-Leffler function and its asymptotic behaviours (stretched exponential for short times and power law for long times).
Presumably that Scott-Blair had guessed the behaviour of the M-L function but he did not have the mathematical background being overall an experimentalist.

Both Scott Blair's model (\ref{Nutting-frac}) and Gerasimov's model (\ref{Gerasimov}) are naturally considered later as special cases of {\it fractional Maxwell's model} with rheological constitutive equation of the form
\begin{equation}
\label{Maxwell-frac}
\sigma(t)  + \lambda^{\alpha} D_{0+, t}^{\alpha} \sigma(t) = E \lambda^{\beta} D_{0+, t}^{\beta} \varepsilon(t),\;\;\; 0 \leq \alpha \leq \beta \leq 1,
\end{equation}
where $E$ is the shear modulus, and $\lambda$ is the relaxation time. This equation generalizes celebrating Maxwell equation in which for the first time
Newtonian law for viscous fluid and Hook's law for elastic solid are combined to describe the behaviour of visco-elastic media
\begin{equation}
\label{Maxwell-classical}
\sigma(t) + \tau D_{t} \sigma(t) = E \tau D_{t} \varepsilon(t).
\end{equation}

Partial case of fractional Maxwell's model is the so-called {\it three-parametric generalized Maxwell's model} with constitutive equation of the type
\begin{equation}
\label{Maxwell-frac3}
\sigma(t)  + a_1 D_{0+, t}^{\alpha} \sigma(t) = b_{0} \varepsilon(t),\;\;\; 0 < \alpha < 1.
\end{equation}
Another popular fractional model with three parameters is the
{\it Kelvin-Voigt fractional model} that presumably for the first time was introduced by Caputo
\cite{Cap67}
in 1967,
\begin{equation}
\label{Voigt-frac}
\sigma(t)  = b_{0}  \varepsilon(t) + b_{1} D_{0+, t}^{\alpha} \varepsilon(t),\;\;\; 0 < \alpha < 1.
\end{equation}
It is a generalization of the classical Kelvin model having the following constitutive equation
\begin{equation}
\label{Kelvin-classical}
\sigma(t)  = E\left[ \sigma(t) + \tau D_{t} \varepsilon(t)\right].
\end{equation}

More general constitutive equation corresponds to the  so called {\it  fractional Zener model}:
\begin{equation}
\label{Zener-frac4}
\sigma(t)  + a_1 D_{0+, t}^{\alpha} \sigma(t) = b_{0}  \varepsilon(t) + b_{1} D_{0+, t}^{\alpha} \varepsilon(t),\;\;\; 0 < \alpha < 1.
\end{equation}
formerly introduced in 1971 by Caputo and Mainardi \cite{CapMai71a}.
Theoretical background for this was done by Bagley and Torvik, see \cite{BagTor83}, \cite{BagTor86}. It has to be pointed out that the above considered Rabotnov's model
(\ref{Rabotnov}) is equivalent  to  the fractional  Zener model, see \cite{Uch08}.

Sometimes the {\it Poynting-Thomson fractional model} is discussed with rheological constitutive equation of the type
\begin{equation}
\label{Zener-frac5}
\sigma(t)  + \frac{E}{E_0} \lambda^{\alpha - \gamma} D_{0+, t}^{\alpha - \gamma} \sigma(t) + \frac{E}{E_0} \lambda^{\beta - \gamma} D_{0+, t}^{\beta - \gamma} \sigma(t) = E \lambda^{\alpha}  D_{0+, t}^{\alpha} \varepsilon(t) + E \lambda^{\beta} D_{0+, t}^{\beta} \varepsilon(t),
\end{equation}
($0 \leq \gamma \leq \alpha \leq \beta \leq 1.$)

More extended discussion of the fractional models in rheology can be found in
\cite{CapMai71b},
\cite{Mai10}, \cite{Pod99}, \cite{SchFriBlu00}, \cite{WBG03}.

\section{The second generation of fractional modelling in rheology}

The development of the fractional differential approach in rheology is associated
with such names as Scott-Blair, Bagley and Torvik, Caputo, Gorenflo and Mainardi \cite{MaGo07},
Friedrich, Schiessel and  Blumen \cite{Fri91}, \cite{SchFriBlu00},  Metzler, Nonnenmacher, Gl\"ockle, Klafter
and Schlesinger \cite{MeGlNo94}, \cite{MetKla04}, Koeller \cite{Koe10}, Podlubny and Heymans \cite{HeyPod06},
Rossikhin and Shitikova \cite{MPPR71}, \cite{RoSh07}, \cite{RoSh13}, and others.

In a series of papers (see, {\it e.g.} \cite{BagTor83}, \cite{BagTor86}) Bagley and Torvik extended the ideas of
Gemant, and of Caputo--Mainardi \cite{CapMai71b}. They have shown that the complex modulus of many materials can be approximated by fractional powers in the frequency. They proposed a general fractional
constitutive relation describing visco-elastic behaviour in its different appearance.

This approach has been successfully
applied to describe rheological behaviour of organic glasses, elastomers,
polyurethane, polyisobutylene, monodisperse polybutadiene and solid amorphous
polymers in a wide temperature range (see for example \cite{AlcMarV99}
and references therein). In \cite{VanUlm06} and \cite{ShaUlm09} have been derived equations
governing the time-dependent indentation response
for axisymmetric indenters into a fractional viscoelastic half-space and have proposed
an original method for the inverse analysis of fractional viscoelastic properties
and applied to experimental indentation creep data of polystyrene. The method
is based on fitting the time-dependent indentation data (in the Laplace domain) to the
fractional viscoelastic model response. It is shown that the particular time-dependent
response of polystyrene is best captured by a bulk-and-deviator fractional viscoelastic
model of the Zener type. We shall dwell in details on fractional differential models
of viscoelasticity and then consider a few standard hydrodynamic problems in
the simplest model of this type.

It is impossible to describe the modern state in the fractional rheology. We refer interested readers to the recent
monographs \cite{Mai10}, \cite{Uch13b}, and to the survey paper \cite{VaMaKi14} for some additional comments on pioneering works in applications of fractional calculus.

\section*{Acknowledgements.} One of the inspiration of this article was a possibility to work at
the G.W.~Scott Blair  reading room at the Aberystwyth University.
The authors are grateful to the official Custodian of the collection,
Prof. S.~Cox, for supporting the
idea of the article and for fruitful discussion of the subject. The work by S.R. is partly supported by
PEOPLE IAPP Project PIAP-GA-2009-251475 HYDROFRAC.



\section*{Appendix}

\noindent {\bf From mathematical curiosity to practical significance}  

\noindent {\em Rheology Bulletin}, {\bf 55}, No. 2, 45--46  (2014)

\vspace{5mm} In March, mathematicians in Aberystwyth were treated to an historical review of
GW Scott Blair's contributions to both philosophy and rheology, with particular
emphasis on his development of fractional calculus from a mathematical curiosity
to a tool of practical significance to rheology. The speaker, Professor Sergei
Rogosin, is a Senior Marie Curie Fellow in the Department of Mathematics and
Physics at Aberystwyth, on leave from Minsk University in Belarus. Professor
Rogosin has worked on fractional calculus previously, and serendipitously
discovered the BSR's Scott Blair Collection on arriving in Aberystwyth. He used
the collection to discover more about Scott Blair's work on fractional calculus, and
his talk was based on a recent paper with Francesco Mainardi (Bologna) \cite{RogMai14}.

As readers may know, not only did Scott Blair help found the BSR, and serve as
one of its first presidents, he also contributed to the founding of the (American)
Society of Rheology twenty years earlier. His interest in fractional calculus was
motivated by trying to explain his experimental results on food rheology, trying to
quantify effects such as "firmness" and "taste" and the influence of material
memory on rheological response \cite{ScoB47b}. Gareth McKinley spoke about the same topic
in Chicheley Hall last year at the INNFM meeting \cite{McKJai} (see also \cite{JaiMcK13}). 
Professor Rogosin explained how Scott Blair picked up on Nutting's theoretical
work, and Gemant's experiments some time later, indicating that fractional
exponents in constitutive relations offered a way to interpolate between Newtonian
liquids and Hookean solids, and consequently introduced the idea of quasiproperties.
The speaker made the link with Volterra integral equations, which had
been developed earlier in the 20th century, providing a robust mathematical
framework in which to explore these relationships (although Scott Blair himself
probably did not recognise this). At the same time as Scott Blair was working in
this area, both Gerasimov and Rabotnov were tackling this problem of
interpolation in similar ways and publishing in the Russian literature (although
Scott Blair spoke Russian, the 1940s were of course a difficult time to exchange
scientific information!). In particular, Rabotnov recognised that using fractional
exponents offered the possibility of predicting power-law decay in time, rather than
exponential decay, allowing much better fitting of the rheological behaviour of
complex fluids. He later acknowledged Scott Blair's contributions in a textbook
published in the 1980s; indeed, Scott Blair's contributions were widely valued both
by experts in fractional modelling \cite{Koe10} and by rheologists \cite{Dor02}.

The seminar was followed by a discussion of nonlinear elasticity and integral
equations, upon which I shall not report here.

The Scott Blair Collection of books on rheology continues to grow thanks to a
financial bequest to the BSR. Donations are also welcome. Further details are
available at https://www.aber.ac.uk/en/is/collections/scottblair/

\vspace{3mm} Simon Cox, Aberystwyth, March 2014
\end{document}